\documentclass[a4paper,11pt,twoside]{article}

\usepackage{xypic}
\usepackage{latexsym}
\usepackage{amsmath}
\usepackage{amstext}
\usepackage{amsfonts}
\usepackage{amssymb}
\usepackage{amsbsy}
\usepackage{eucal}
\usepackage{enumerate}
\usepackage{theorem}

\xyoption{all}
\CompileMatrices

\newtheorem{lemma}{Lemma}
\newtheorem{proposition}{Proposition}

\theorembodyfont{\rmfamily}
\newtheorem{remark}{Remark}

\newcommand{\id}{\mathop{\rm id}\nolimits}

\newcommand{\Spec}{\mathop{\rm Spec}\nolimits}
\newcommand{\End}{\mathop{\rm End}\nolimits}

\newcommand{\largefield}{K}

\newcommand{\CPmap}{\phi}

\newcommand{\qed} {\hbox{} \nolinebreak \hfill $\;\Box$}
\newcommand{\dotcup}{\mathrel{\Dot \cup}}

\begin{document}

\addtolength{\baselineskip}{0.5 pt}

\setlength{\parskip}{0.2 ex}

\title{Ordinary elliptic curves of high rank over $\overline{\mathbb{F}}_p(x)$ with constant $j$-invariant II}

\author{Claus Diem and Jasper Scholten}

\date{\today}

\maketitle

\begin{abstract}

We show that for all odd primes $p$, there exist ordinary elliptic
curves over $\overline{\mathbb{F}}_p(x)$ with arbitrarily high rank and
constant $j$-invariant. This shows in particular that there are elliptic curves with arbitrarily high rank over these fields for which the corresponding elliptic surface is not
supersingular. The result follows from a theorem which states
that for all odd prime numbers $p$ and $\ell$, there exists a
hyperelliptic curve over $\overline{\mathbb{F}}_p$ of genus
$(\ell-1)/2$ whose Jacobian is isogenous to the power of one ordinary
elliptic curve.
\end{abstract}

\vspace{1 ex}
\begin{center}\begin{minipage}{110mm}\footnotesize{\bf Key words:} Elliptic curves of high rank, Jacobians
\end{minipage}
\end{center}

\begin{center}\begin{minipage}{110mm}\footnotesize{\bf MSC2000:} Primary: 11G05; Secondary: 11G20, 14H40, 14H52.

\end{minipage}
\end{center}
\vspace{1 ex}

\pagestyle{myheadings} \markboth{\sc Diem, Scholten}{\sc Ordinary elliptic curves of high rank II}

\section{Introduction}\label{intro}

The purpose of this work is to prove the following theorem and its corollary.

\paragraph{Theorem}
\emph{Let $p$ and $\ell$ be odd prime numbers. Then there exists an ordinary hyperelliptic curve $H$ over $\overline{\mathbb{F}}_p$ of genus $(\ell -1)/2$ whose Jacobian variety is isogenous to the power of one ordinary elliptic curve.}

\paragraph{Corollary}
\emph{Let $p$ be a fixed odd prime number. Then there exist ordinary elliptic curves over $\overline{\mathbb{F}}_p(x)$ with arbitrarily high rank and constant $j$-invariant.}

\medskip

The corollary follows from the theorem and the following general fact which was first pointed out by Shafarevich and Tate (\cite{shaftate}): Let $H$ be a hyperelliptic curve over some field $k$ with a $k$-rational divisor of degree 1, and let $E$ be an elliptic curve over $k$ such that $J_H \sim E^r \times A$ for some abelian variety $A$ over $k$. Let $k(H)|k(x)$ be an extension of degree 2. Then the quadratic twist of $E_{k(x)}$ with respect to $k(H)|k(x)$ has rank at least $r$. For details we refer to the introduction of \cite{BDS}.

Using the construction in \cite{GHS}, we have shown in \cite{BDS} that for every $r \in \mathbb{N}$ there exists an ordinary hyperelliptic curve over $\mathbb{F}_{2^r}$ whose Jacobian variety is completely decomposable and contains up to isogeny a factor $E^r$ for some ordinary elliptic curve $E$ (\cite[Theorem 1]{BDS}). This implies that the corollary holds in fact for all primes $p$. We note that in \cite{BDS} the corollary 
has already been proven under the assumption of Artin's Primitive Root Conjecture.

\medskip

It is not known whether there exist non-constant elliptic curves over $\mathbb{C}(x)$ with arbitrarily high rank. On the other hand, as noted in \cite{BDS}, there are several known classes of non-constant
elliptic curves with arbitrarily high rank over rational function fields
in positive characteristic. Previous
to our work \cite{BDS}, examples of both high-rank
supersingular elliptic curves and high-rank 
elliptic curves with non-constant $j$-invariant (which are thus ordinary) have
been given. 
But all these examples are related to supersingular
varieties because 
the corresponding elliptic surfaces are
supersingular (in the sense of Shioda, that is, the $\ell$-adic 
second cohomology groups
are
generated by divisor classes, cf.\ \cite[page 235]{shioda2}).
This is somewhat
of a weakness if one wants to view the positive characteristic examples as
an indication of what should be true in characteristic~0.

Just as the elliptic curves of high rank in \cite{BDS}, the elliptic curves of high rank given in this work do not correspond to supersingular surfaces.

\medskip

In \cite{ES} it was asked by Ekedahl and Serre whether there exist curves in characteristic zero of arbitrarily high genus with completely decomposable Jacobian. It is well-known that the corresponding question for any fixed positive characteristic can be answered affirmatively: In every positive characteristic there exist supersingular curves of arbitrarily high genus, and over an algebraically closed field, the Jacobians of these curves are completely decomposable. Together with \cite[Theorem 1]{BDS} our theorem shows that in every positive characteristic there also exist ordinary curves of arbitrarily high genus with completely decomposable Jacobian.

\subsubsection*{Outline of the proof of the theorem}

We consider families of hyperelliptic curves which have already been studied by Brumer and Mestre (see \cite{mestre} and reference [3] in \cite{mestre}). Over an algebraically closed field, the Jacobian varieties of the members of these families have real multiplication. Via a degeneration argument we give a criterion under which the families are generically ordinary. If we consider certain members of these families such that the real multiplication is not defined over the ground field, we can conclude that the Jacobian varieties are isogenous to Weil restrictions of ordinary elliptic curves (here we use \cite[Proposition 2]{BDS}). The existence of suitable curves is guaranteed by a density argument.

\subsubsection*{Terminology and notations}
Throughout the work, we use the following terminology concerning curves: If not stated otherwise, a \emph{curve} over a field is assumed to be proper and geometrically reduced but not necessarily smooth or irreducible. A (relative) curve over a locally noetherian connected base $S$ is a proper and flat morphism $C \longrightarrow S$ whose fibers are curves. A \emph{semi-stable} curve over $S$ is a curve over $S$ such that the only singularities of the fibers are ordinary double points.

We denote a relative curve $C$ over $S$ by $C/S$. For $s \in S$ or a geometric point $s$ of $S$, we denote the (geometric) fiber of $C$ over $s$ by~$C_s$. If $C/S$ and $D/S$ are curves, we call a finite and flat $S$-morphism $C \longrightarrow D$ a \emph{cover}.

Note that the Euler-Poincar\'e characteristic of the fibers of a relative curve is constant (see the first corollary to the Theorem in \cite[Chapter II, \S 5]{Mu}), and thus so is the (arithmetic) genus; we define the \emph{genus} of a relative curve to be the genus of any of its fibers.

If $A$ is an abelian variety over a field $k$, we denote the ring of endomorphisms of $A$ \emph{over} $k$ by $\End(A)$.

\section{The curves}
\label{curves-alg-closed-fields}
In this section, we define the hyperelliptic curves in question as well as the the real multiplication on the Jacobians of these curves.

Let $\ell > 1$ be an odd prime. Let us fix a ground field $k$ of any
characteristic $\neq 2$ with algebraic closure $\overline{k}$. All
curves in this section are over the fixed field $k$ (but the real
multiplication we want to define might not be defined over $k$). 
We assume that we are given an \'etale isogeny of elliptic curves 
$f : \tilde{E} \longrightarrow E$ over $k$ of degree $\ell$. Let $\rho : E \longrightarrow \mathbb{P}^1_k, \tilde{\rho} : \tilde{E} \longrightarrow \mathbb{P}^1_k$ be covers of degree 2 which are ramified at the 2-torsion points. As the covers $\rho, \tilde{\rho}$ are in fact geometric quotients of $E$ and $\tilde{E}$ respectively by the automorphisms $[-1]$ on these curves, we have an induced cover $u : \mathbb{P}^1_k \longrightarrow \mathbb{P}^1_k$ of degree $\ell$ with $\rho \circ f = u \circ \tilde{\rho}$.

Let $P := \rho(0_E)$, where $0_E$ is the zero point of $E$. Let $a : \mathbb{P}^1_k \longrightarrow \mathbb{P}^1_k$ be a cover of degree 2 which is ramified at $P$ and some other point $P' \notin \rho(E[2])$. 

Now let $D$ be the normalization of the fiber product of $a : \mathbb{P}^1_k \longrightarrow \mathbb{P}^1_k$ and $u : \mathbb{P}^1_k \longrightarrow \mathbb{P}^1_k$, and let $d : D \longrightarrow \mathbb{P}^1_k, \; \tilde{a} : D \longrightarrow \mathbb{P}^1_k$ be the covers which are induced by the projections to the first and second factor respectively. We thus have the following commutative diagram.
\begin{equation}
\label{D-def}
\xymatrix{
D \ar^{d}[dd] \ar_(0.3){\tilde{a}}[drr] & {\tilde{E}} \ar^(0.6){f}|(0.25)\hole[dd] \ar^{\tilde{\rho}}[dr] \\
& & {\mathbb{P}^1_k} \ar^{u}[dd] \\
{\mathbb{P}^1_k} \ar_(0.3){a}[drr] & E \ar^{\rho}[dr] \\
& & {\mathbb{P}^1_k}}
\end{equation}
Obviously $D$ is hyperelliptic. Let us calculate its genus: 

By Abhyankar's Lemma (cf.\ \cite[Proposition III.8.9]{Stich}), the cover $\tilde{a}_{\overline{k}} : D_{\overline{k}} \longrightarrow \mathbb{P}^1_{\overline{k}}$ is unramified outside the preimages of $P$ and $P'$ under $u_{\overline{k}}$. Over the point $P$ there is the point $\rho'(0_{\tilde{E}})$ which is unramified and $(\ell-1)/2$ other points in $\mathbb{P}^1(\overline{k})$ which are ramified of order 2. The point $P'$ is unramified with respect to $u_{\overline{k}}$.

Again by Abhyankar's Lemma, the cover $\tilde{a}_{\overline{k}} : D_{\overline{k}} \longrightarrow \mathbb{P}^1_{\overline{k}}$ is ramified exactly at $\tilde{\rho}(0_{\tilde{E}})$ and at the preimages of $P'$. This means that there are exactly $\ell +1$ ramification points, and the genus of $D$ is $(\ell-1)/2$.

\medskip

To introduce the real multiplication, we first consider the normalization $C$ of $\mathbb{P}^1_k \times_{\mathbb{P}^1_k} E$, where the product is taken with respect to $a$ and $\rho$. Note that $C$ is a genus 2 curve. Let $\tilde{C} := C \times_E \tilde{E}$. Note that as $f : \tilde{E} \longrightarrow E$ is an \'etale cover of degree $\ell$, so is $\tilde{C} \longrightarrow C$. In particular, $\tilde{C}$ is smooth, that is, it is a smooth and irreducible curve (of genus $\ell +1$). We have the following commutative diagram in which we also introduce some names for the covers.
\begin{equation}
\label{C-tilde-def}
\xymatrix{
{\tilde{C}} \ar_{c}[dd] \ar_{\tilde{\CPmap}}[dr] \ar[drr] \\
& D \ar^{d}[dd] \ar_(0.3){\tilde{a}}[drr] & {\tilde{E}} \ar^(0.6){f}|(0.25)\hole[dd] \ar^{\tilde{\rho}}[dr] \\
C \ar_{\CPmap}[dr] \ar|\hole[drr] & & & {\mathbb{P}^1_k} \ar^{u}[dd] \\
& {\mathbb{P}^1_k} \ar_(0.3){a}[drr] & E \ar^{\rho}[dr] \\
& & & {\mathbb{P}^1_k}}
\end{equation}
Let $\largefield|k$ be a field extension such that $\tilde{E}_{\largefield}$ has a non-trivial translation $\tau_{\tilde{E}_{\largefield}}$ (of order $\ell$) which is an automorphism of the cover $f_{\largefield} : \tilde{E}_{\largefield} \longrightarrow E_{\largefield}$ induced by $f$. Then $\tau_{\tilde{E}_{\largefield}}$ induces by base-change an automorphism of order $\ell$ of $\tilde{C}_\largefield = (C \times_E \tilde{E})_{\largefield}$. Note that $\tau_{\tilde{C}_{\largefield}}$ generates the automorphism group of the cover $c_{\largefield}: \tilde{C}_{\largefield} \longrightarrow C_{\largefield}$.

Let $A$ be the reduced connected component of the zero of $\ker(c_*) < J_{\tilde{C}}$. Note that $\id + \tau_{\tilde{C}_{\largefield}}^* + \cdots + \tau_{\tilde{C}_{\largefield}}^{*\ell-1}$ is defined over $k$, and $A$ is the reduced connected component of the zero of $\ker(c^* c_*) = \ker(\id + \tau_{\tilde{C}_\largefield}^* + \cdots + \tau_{\tilde{C}_\largefield}^{*\ell-1})$. Obviously $\tau^*_{\tilde{C}_\largefield}(A_\largefield) \subseteq A_\largefield$ and the minimal polynomial of $\tau_{\tilde{C}_\largefield}^*|_{A_\largefield}$ divides $X^{\ell-1} + \cdots + X + 1$. As the latter polynomial is irreducible in $\mathbb{Q}[X]$ and $\tau^*_{\tilde{C}_\largefield}|_{A_\largefield}$ operates non-trivially on $A_\largefield$, we conclude:

\begin{lemma}
\label{A-decomp}
The ring $\mathbb{Z}[\tau_{\tilde{C}_\largefield}^*|_{A_\largefield}]< \End(A_\largefield)$ is isomorphic to $\mathbb{Z}[\zeta_\ell]$ (with $\tau _{\tilde{C}_\largefield}^*|_{A_\largefield} \longleftrightarrow \zeta_\ell$).
\end{lemma}

\begin{lemma}
$\tilde{\CPmap}^*(J_D)$ is contained in $A$.
\label{J_H<A}
\end{lemma}
\emph{Proof.} We proof this after base-change to $\largefield$. Let $\delta$ be the non-trivial automorphism of $\tilde{C}$ over $D$ (of order 2). Note that $\tau_{\tilde{C}_{\largefield}}$ and $\delta_{\largefield}$ generate the automorphism group of $\CPmap_{\largefield} \circ c_{\largefield} : \tilde{C}_{\largefield} \longrightarrow \mathbb{P}^1_{\largefield}$ which is in fact a dihedral group. In particular, we have
\begin{equation}
\label{add-everything->trivial}
\begin{array}{c}
(\id + \tau_{\tilde{C}_{\largefield}}^* + \cdots  + \tau_{\tilde{C}_{\largefield}}^{*\ell-1}) \circ (\id + \delta_{\largefield}^*)  = \\
\id + \tau_{\tilde{C}_{\largefield}}^* + \cdots  + \tau_{\tilde{C}_{\largefield}}^{*\ell-1} + \delta_{\largefield}^* + \tau_{\tilde{C}_{\largefield}}^* \circ  \delta_{\largefield}^* + \cdots  + \tau_{\tilde{C}_{\largefield}}^{*\ell-1} \circ \delta_{\largefield}^*  = \\
c_{\largefield}^* \circ \CPmap_{\largefield}^* \circ \CPmap_{\largefield*} \circ c_{\largefield*} = 0 \, .
\end{array}
\end{equation}

Now $\tilde{\CPmap}_{\largefield}^*(J_{D_{\largefield}}) =  \tilde{\CPmap}^*_{\largefield} \tilde{\CPmap}_{\largefield*} (J_{\tilde{C}_{\largefield}}) = (\id + \delta_{\largefield}^*)(J_{\tilde{C}_{\largefield}})$. By (\ref{add-everything->trivial}), this lies in $\ker(\id + \tau_{\tilde{C}_{\largefield}}^* + \cdots + \tau_{\tilde{C}_{\largefield}}^{*\ell-1})$. \qed
\smallskip

Now let $\tau := \tilde{\CPmap}_{\largefield*} \circ \tau_{\tilde{C}_{\largefield}}^* \circ \tilde{\CPmap}_{\largefield}^* \in \End(J_{D_{\largefield}})$.

\begin{proposition}
\label{real-mult}
$\mathbb{Z}[\tau] < \End(J_{D_{\largefield}})$ is isomorphic to $\mathbb{Z}[(\tau_{\tilde{C}_{\largefield}}^* + \tau_{\tilde{C}_{\largefield}}^{*-1})|_{A_{\largefield}}] < \End(A_{\largefield})$ as rings and Galois modules (with $\tau \longleftrightarrow (\tau_{\tilde{C}_{\largefield}}^* + \tau_{\tilde{C}_{\largefield}}^{*-1})|_{A_{\largefield}}$). In particular, the ring $\mathbb{Z}[\tau]$ is isomorphic to $\mathbb{Z}[\zeta_\ell + \zeta_\ell^{-1}]$.
\end{proposition}
\emph{Proof.} This follows by an analogous argument to the one presented on the lower part of page 498 of \cite{BDS}.
\qed

\begin{remark}
For the proof of the theorem, we do not need the real multiplication on $J_D$ but merely the automorphism $\tau^*_{\tilde{C}_{\largefield}}|_{A_{\largefield}}$. This is similar to the situation in \cite{BDS}.
\end{remark}

\begin{remark}
\label{comp-curve}
The automorphism group of $a \circ \CPmap : C \longrightarrow \mathbb{P}^1_k$ is $(\mathbb{Z}/2 \mathbb{Z})^2$, thus besides the quotients $C \longrightarrow E$ and $\CPmap : C \longrightarrow \mathbb{P}^1_k$ of degree 2, there is another quotient $C \longrightarrow E'$ of degree 2 together with an induced cover $\rho' : E' \longrightarrow \mathbb{P}^1_k$ of degree 2. Now, $\rho' : E' \longrightarrow \mathbb{P}^1_k$ is ramified at $\rho(E[2] - \{ 0_E \} ) \cup \{ P' \}$, thus $E'$ has genus 1. By fixing the preimage of $P'$ as zero-point, $E'$ is turned into an elliptic curve. Obviously $C$ is also the normalization of $E \times_{\mathbb{P}^1_k} E'$. In the next section this observation serves as the starting point.
\end{remark}

\section{Generically ordinary families}
In this section, we introduce families of curves (i.e.,\ relative curves) $\tilde{C}/S$ whose fibers are -- except for some ``degenerate'' fibers -- curves denoted by $\tilde{C}$ in the previous section. We strive for proving that under some conditions, these families are generically ordinary, that is, their generic fibers are ordinary. Note that as ``ordinariness'' is an ``open property'' (cf.\ \cite[Proposition 1.5]{Bouw-Luminy}), this means that there exists a non-empty open subset $U \subset S$ such that the fibers over $U$ are ordinary. To prove this, we force that some fibers are not irreducible.

The setting in this section is more general than what is actually needed for the proof of the theorem.
\smallskip

Let $S$ be an integral, regular scheme over $\mathbb{Z}[1/2]$ of dimension 1, and let $E/S$ and $E'/S$ be two (relative) elliptic curves. Let $[0_E]$ be the Cartier divisor associated to the zero-section, and let $E[2]^{\#} := E[2] - [0_E]$. (Similar definitions for $E'$.)

We assume that we are furthermore given
\begin{itemize}
\item
an isomorphism $\psi : E[2] \tilde{\longrightarrow} E'[2]$ of group schemes over $S$,
\item
an \'etale isogeny of elliptic curves $f : \tilde{E} \longrightarrow E$ of some degree $\ell$ (which need neither be odd nor prime).
\end{itemize}

Let $\rho : E \longrightarrow E/[-1]$, $\rho' : E' \longrightarrow E'/[-1]$ be the canonical maps. Note that $E/[-1]$ and $E'/[-1]$ are both $\mathbb{P}^1$-bundles (in fact if $q : E \longrightarrow S$ is the structure map then $E/[-1] \simeq \mathbb{P}(q_* \, \mathcal{L}(2[0_E]))$).  Furthermore, both $E[2]^{\#}$ and $E'[2]^{\#}$ are \'etale over $S$ of degree 3 and $\rho|_{E[2]}$ as well as $\rho'|_{E'[2]}$ are closed immersions. By \cite[Proposition B.4]{Di-degree2} there exists a unique $S$-isomorphism
\[ \gamma : E/[-1] \tilde{\longrightarrow} E'/[-1] \]
with $\gamma \circ \rho|_{E[2]^{\#}} = \rho' \circ \psi|_{E[2]^{\#}}$.

\smallskip

We identify $E/[-1]$ and $E'/[-1]$ via $\gamma$ and let $\mathbf{P} := E/[-1] = E'/[-1]$, $P := \rho([0_{E}])$ and $P' := \rho'([0_{E'}])$. 

\smallskip

We have the following easy lemma (see \cite[Lemma A.1]{Di-degree2}).

\begin{lemma}
Let $s$ be a geometric point of $S$. Then the following are equivalent:
\begin{itemize}
\item
There exists an isomorphism $\alpha : E_s \longrightarrow E'_s$ with $\alpha|_{E_s[2]} = \psi_s$.
\item
$P_s = P'_s$.
\end{itemize}
\end{lemma}

We now make the following assumption:

\paragraph{Assumption}
The divisors $P$ and $P'$ are not equal (i.e.,\ they are not equal at the generic point of $S$) but there exists some point $s \in S$ such that $P_s = P'_s$ and $E_s (=E'_s)$ is ordinary.

\begin{remark}
For the proof of the theorem 
we will deal with the special case where $S$ is an affine curve over $\mathbb{F}_q$ and $E'/S$ is a constant family.
\end{remark}

\begin{lemma}
$E \times_{\mathbf{P}} E'$ is integral, and $\kappa(E \times_{\mathbf{P}} E') \simeq \kappa(E) \otimes_{\kappa(\mathbf{P})} \kappa(E')$.
\end{lemma}
\emph{Proof.} 
The ring $E \otimes_{\kappa(\mathbf{P})} \kappa(E')$ is a field because by assumption, the generic points of $\rho'([0_{E'}])$ and $\rho([0_E])$ are distinct.

Let $A$ be the coordinate ring of an affine open part $U$ of $\mathbf{P}$, let $B$ and $B'$ be the coordinate rings of the preimages of $U$ in $E$ and $E'$. We claim that the canonical map $B \otimes_A B' \longrightarrow \kappa(B) \otimes_{\kappa(A)} \kappa(B') \simeq \kappa(E) \otimes_{\kappa(\mathbf{P})} \kappa(E')$ is injective.

We have $\kappa(B) \otimes_{\kappa(A)} \kappa(B') \simeq (B \otimes_A \kappa(A)) \otimes_{\kappa(A)} (B' \otimes_A \kappa(A)) \simeq (B \otimes_A B') \otimes_A \kappa(A)$ as $B$ and $B'$ are finite over $A$. We thus have to show that the map $B \otimes_A B' \longrightarrow (B \otimes_A B') \otimes_A \kappa(A)$ is injective. Now, $A \longrightarrow \kappa(A)$ is injective and $B \otimes_A B'$ is flat over $A$. This implies that $B \otimes_A B' \longrightarrow (B \otimes_A B') \otimes_A \kappa(A)$ is injective. It follows that $B \otimes_A B'$ is reduced and that $\kappa(B) \otimes_{\kappa(A)} \kappa(B') \simeq \kappa(E) \otimes_{\kappa(\mathbf{P})} \kappa(E')$ is its function field.
\qed


\medskip

Let $C$ be the normalization of $E \times_{\mathbf{P}} E'$, and let $\pi : C \longrightarrow E$ and $\pi' : C \longrightarrow E'$ be the maps which are induced by the projections to the factors.

\begin{lemma}
The morphisms $\pi$ and $\pi'$ as well as the normalization morphism $C \longrightarrow E \times_{\mathbf{P}} E'$ are finite, in particular $C/S$ is proper.
\end{lemma}
\emph{Proof.}
We first note that as $S$ is regular and $E$ is smooth over $S$, $E$ is also regular and in particular normal (see \cite[Expos\'e II, Proposition 3.1]{SGA}, \cite[Theorem 19.4]{Ma}).

The projection $E \times_{\mathbf{P}} E' \longrightarrow E$ is finite because it is induced by base-change from the finite morphism $E' \longrightarrow \mathbf{P}$. This implies in particular that the normalization of $E \times_{\mathbf{P}} E'$ is equal to the normalization of $E$ in $\kappa(E) \otimes_{\kappa(\mathbf{P})} \kappa(E')$. As $E$ is normal and the extension $\kappa(E) \otimes_{\kappa(\mathbf{P})} \kappa(E')|\kappa(E)$ is separable, this implies that $\pi : C \longrightarrow E$ is finite (see \cite[Chapter 4, Proposition 1.25]{Liu}). As $E \times_{\mathbf{P}} E'$ is finite over $E$, implies that the normalization morphism $C \longrightarrow E \times_{\mathbf{P}} E'$ is finite too. Clearly, $\pi'$ is also finite.
\qed

\begin{proposition}
\label{semi-stable}
$C/S$ is a semi-stable genus 2 curve and $\pi, \pi'$ are covers of degree~2.

For $s \in S$ with $P_s \neq P'_{s}$, $C_s$ is an irreducible smooth genus 2 curve. For $s \in S$ with $P_s = P'_s$, $C_s$ is isomorphic to the union of two copies of $E_s (=E'_s)$ intersecting at $0_{E_s}$. One can choose an isomorphism such that $\pi_s$ restricts to $\id_{E_s}$ on both components and $\pi'_s$ restricts to $\id_{E_s}$ on one component and to $-\id_{E_s}$ on the other.
\end{proposition}
\emph{Proof.}
The field $\kappa(S)$ is algebraically closed in $\kappa(E) \otimes_{\kappa(\mathbf{P})} \kappa(E')$ and $S$ is normal. By Grothendieck's results on ``Zariski's Main Theorem'', this implies that the geometric fibers of $C$ over $S$ are connected (see \cite[(4.3.12)]{EGAIII2}).

Let $W$ be the ramification locus of $E \times_{\mathbf{P}} E' \longrightarrow \mathbf{P}$. Then $(E \times_{\mathbf{P}} E') - W$ is normal, because the domain of an \'etale morphism mapping to a normal scheme is normal (see \cite[Expos\'e I, Corollaire 9.11]{SGA}). It follows that $C \longrightarrow E \times_{\mathbf{P}} E'$ induces an isomorphism between the complement of the preimage of $W$ in $C$ and $(E \times_{\mathbf{P}} E') - W$. 

Let $s \in S$. Then $W_s$ is zero-dimensional, and as the canonical morphism $C_s \longrightarrow E_s \times_{\mathbf{P}_s} E'_s$ is finite (by the preceeding lemma), the preimage of $W_s$ in $C$ is zero-dimensial too. Moreover, the irreducible components of $E_s \times_{\mathbf{P}_s} E'_s$ are 1-dimensional (there are one or two irreducible components depending on whether $P_s \neq P_s'$ or $P_s = P_s'$). As $C_s$ is connected, this implies that $C_s \longrightarrow E_s \times_{\mathbf{P}_s} E'_s$ is birational on the irreducible components of $C_s$. The fact that $C/S$ is proper (see preceding lemma) now implies that $C_s$ is proper too. This means that $C_s$ is a curve.

As $S$ is integral of dimension 1 and $C$ is integral, $C/S$ is flat (see \cite[Proposition 9.7]{Ha}). It follows that $C/S$ is a curve.

\smallskip

By Abhyankar's Lemma (\cite[Expos\'{e} X, Lemme 3.6]{SGA}) and ``purity of the branch locus'' (\cite[Expos\'{e} X, Th\'eor\`eme 3.1]{SGA}), $\pi$ is \'etale outside $(\pi')^{-1}([0_{E'}]) = (\rho' \circ \pi')^{-1}(P'_s)$ and $\pi'$ is \'etale outside $\pi^{-1}([0_{E}]) = (\rho \circ \pi)^{-1}(P_s)$. This implies that $C/S$ is smooth outside the preimage of $P \cap P'$.

This means in particular that for $s \in S$ with $P_s \neq P_s'$, the fiber $C_s$ is a non-singular curve. From the assumption it follows now that the generic fiber is a smooth curve. Clearly it has genus 2. This means that $C/S$ is a curve of genus 2.

\smallskip

Now let $s \in S$ with $P_s = P'_s$. 

Note that the normalization of $E_s \times_{\mathbf{P}_s} E'_s = E_s \times_{\mathbf{P}_s} E_s$ is isomorphic to $E_s \dotcup E_s$. The identification can be made such that the first projection map restricts to $\id_{E_s}$ on both components and the second projection map restricts to $\id_{E_s}$ on one component and to $-\id_{E_s}$ on the other. We thus have a canonical map $E_s \dotcup E_s \longrightarrow C_s$ such that the restrictions to each of the two components are closed immersions. 

As outside the preimage of $P_s$ the curve $C_s$ is smooth, we furthermore have an isomorphism outside the preimage of $P_s$. Now, the preimage of $P_s$ in $E_s \times_{\mathbf{P}_s} E'_s$ is $\{ (0_{E_s}, 0_{E'_{s}}) \}$, and the preimage of this in $E_s \dotcup E_s$ is the set consisting of the two zero points.

As $C_s$ is connected (see beginning of the proof), this implies that set-theoretically $C_s$ is obtained from $E_s \dotcup E_s$ by identifying the two zero points, and that moreover, outside the intersection point $C_s$ is smooth. We claim that the intersection point is an ordinary double point.

Note that the genus of $C_s$ is 2 because the genus of the fibers of $C/S$ is constant (see the first corollary to the Theorem in \cite[Chapter II, \S 5]{Mu}). Together with \cite[Chapter 7, Proposition 5.4]{Liu} this implies that the intersection point is an ordinary double point.

This implies that $C_s$ is semi-stable, and hence $C/S$ is semi-stable too.

\smallskip

Clearly the morphisms $\pi$ and $\pi'$ induce covers of degree 2 on the fibers. They are flat by \cite[(11.3.11)]{EGAIV3}, and we have already seen in the preceding lemma that they are finite. They are thus covers of curves of degree 2.
\qed

\medskip

Let $\tilde{C} := C \times_E \tilde{E}$. Note that as $f$ is \'etale, so is the induced map $c : \tilde{C} \longrightarrow C$. It follows that $\tilde{C}/S$ is a semi-stable curve of genus $\ell+1$.

\medskip

The main result of this section is:

\begin{proposition}
\label{gen-ordinary}
The curve $\tilde{C}/S$ is generically ordinary.
\end{proposition}
\emph{Proof.} We assumed that there exists some point $s \in S$ with $P_s = P'_s$ such that $E_s$ is ordinary. Let us fix such a point. Then the irreducible components of $C_s$ are elliptic curves isomorphic to $E_s$. In fact, by the previous proposition, the preimage of such an irreducible component is isomorphic to the ordinary elliptic curve $\tilde{E}_s$. This implies that $\tilde{C}_s$ is ordinary (cf.\ the text below Lemma 1.3 in \cite{Bouw-Luminy}), thus the curve $\tilde{C}/S$ is generically ordinary (cf.\ \cite[Proposition 1.5]{Bouw-Luminy}).
\qed

\begin{remark}
Let $\Omega$ be an algebraically closed field, and let $s \in S(\Omega)$ be a geometric point of $S$ with $P_s = P'_s$. Let $Q_1, \ldots, Q_\ell \in \tilde{E}_s(\Omega)$ be the preimages of $0_{E_s}$ under $f_s$. Then $\tilde{C}_s$ is in fact the union of two copies of $\tilde{E}_s$ where any point $Q_i$ on one component is identified with the corresponding point on the other component. This means that the dual graph of $\tilde{C}_s$ (cf.\ \cite[Section 1]{Bouw-Luminy}) consists of two vertices which are connected by $\ell$ edges. In particular it has $\ell-1$ loops. This is consistent with the genus being $\ell +1$.
\end{remark}

\begin{remark}
\label{D/S}
By \cite[Lemma 5.6]{LK}, the covers $\pi, \pi'$ both have automorphisms of order 2. Let $\sigma_C$ be the composition of these two automorphisms. Then $\sigma_C$ is also an involution, and we have an induced cover $a : C/\langle \sigma_{C} \rangle \longrightarrow \mathbf{P}$ of degree 2 whose branch locus is $P \cup P'$ (see \cite{Di-degree2} for more information). If $s \in S$ such that $P_s \neq P'_s$, then $C/\langle \sigma_{C} \rangle \approx \mathbb{P}^1_{\kappa(s)}$. Note that this corresponds to the setup of the previous section.

Now let $\ell$ be odd. Let $\tilde{\mathbf{P}} := \tilde{E} / [-1]$. Then as in the previous section we have an induced cover $u : \tilde{\mathbf{P}} \longrightarrow \mathbf{P}$ of degree $\ell$. Let $D$ be the normalization of the fiber product $C / \langle \sigma_{S} \rangle \times_{\mathbf{P}} \tilde{\mathbf{P}}$ (which is integral), and let $\tilde{a} : D \longrightarrow \tilde{\mathbf{P}}$ be the canonical map. Let $U$ be the open subscheme of $S$ where $P$ and $P'$ do not meet. Then $D_U/U$ is a hyperelliptic curve of genus $(\ell-1)/2$ in the sense of \cite{LK}.

For this, one has to show that $D_U/U$ is smooth and that $\tilde{a}_U : D_U \longrightarrow \tilde{\mathbf{P}}_U$ is a degree 2 cover. The smoothness follows with Abhyankar's Lemma and ``purity of the branch locus'' similarly to the proof of Proposition \ref{semi-stable}, and the statement on $\tilde{a}_U$ is also proved analogously to the corresponding statement in Proposition \ref{semi-stable}.

Note that for $s \in U$ and prime $\ell$, the curve $D_s$ is a hyperelliptic curve called $D$ in the previous section, i.e.,\ the point of study of this work. We do not need the family $D/S$ in the sequel, though.
\end{remark}

\section{Modular curves and proof of the theorem}

Let $p$ and $\ell$ be odd primes.

We want to apply \cite[Proposition 2]{BDS} to the abelian variety $A < J_{\tilde{C}}$ for some curve $\tilde{C}$ as in Section \ref{curves-alg-closed-fields} over some finite field of characteristic $p$.

For this, we first fix an ordinary elliptic curve $E'$ over a finite field $\mathbb{F}_q$ of characteristic $p$. By enlarging the base field if necessary, we assume that the two-torsion points are all rational, and that $\zeta_4 \in \mathbb{F}_q$. We fix a cover $\rho' : E' \longrightarrow \mathbb{P}^1_{\mathbb{F}_q}$ of degree 2, ramified as always at the 2-torsion points.

Note that the moduli problem ``\'etale isogenies of order $\ell$ of elliptic curves with rational $4$-torsion and a fixed determinant over $\mathbb{F}_q$-schemes'' is (finely) represented by an irreducible affine curve over $\mathbb{F}_q$. Let us call this curve $Y_0$. Let $Y_1$ be the affine modular curve parameterizing elliptic curves with rational 4-torsion and a fixed determinant with an automorphism of order $\ell$ over $\mathbb{F}_q$-schemes. Note that we have a canonical \'etale Galois cover $Y_1 \longrightarrow Y_0$ with Galois group $(\mathbb{Z}/\ell\mathbb{Z})^*$. Thus in particular the cover is cyclic of order $\ell-1$.

Let $f: \tilde{E} \longrightarrow E$ be a universal \'etale isogeny of order $\ell$ over $Y_0$. The curve $Y_0$ will play the role of the base-scheme, denoted $S$ in the previous section. The isogeny $f$ will play the same role as in the previous section. The complementary elliptic curve for the construction is $E'_{Y_0}$, obtained by base-change from $\Spec(\mathbb{F}_q)$ to $Y_0$.

By definition of $E$ and $E'$, both $E[2]$ and $E'_{Y_0}[2]$ are isomorphic to the group scheme $(\mathbb{Z}/2\mathbb{Z})^2$ over $Y_0$. There thus exists an isomorphism between $E[2]$ and $E'_{Y_0}[2]$ as $S$-group schemes. Let us choose such an isomorphism. We are then in the situation of the previous section, and we define the curves $C/Y_0$ and $\tilde{C}/Y_0$ as in the previous section.

The ``Assumption'' of the previous section is obviously satisfied, and by Proposition \ref{gen-ordinary}, $\tilde{C}/Y_0$ is generically ordinary. Recall that this means that there are only finitely many geometric fibers which are not ordinary.

Now consider the cover $Y_1 \longrightarrow Y_0$ of degree $\ell-1$. Note that by definition of $Y_1$, if $y \in Y_0$, $y_1$ is a preimage of $y$ and $k := \mathbb{F}_q(y)$, $\largefield := \mathbb{F}_q(y_1)$ are the corresponding residue fields, there exists a non-trivial translation $\tau_{\tilde{E}_\largefield}$ (of order $\ell$) on $(\tilde{E}_y)_{\largefield}$ which is an automorphism of the cover $(f_y)_\largefield : (\tilde{E}_y)_\largefield \longrightarrow (E_y)_\largefield$, and $\largefield$ is the smallest extension of $k$ over which $\tau_{(\tilde{E}_y)_\largefield}$ is defined.

By the Chebotarev Density Theorem (see e.g.\ \cite{Jarden-Chebotarev}), there exist infinitely many $y \in Y_0$ which are inert. In particular, there exists such a $y$ such that $P_y \neq P'_y$ (which implies that $C_y$ and $\tilde{C}_y$ are smooth and irreducible) and $\tilde{C}_y$ is ordinary. Let us fix such a $y$. Let $y_1$ be the unique preimage of $y$ in $Y_1$. Let $k := \mathbb{F}_q(y)$ and $\largefield := \mathbb{F}_q(y_1)$. Note that by the choice of $y$, $[\largefield : k] = \ell-1$.

Now similarly to Section \ref{curves-alg-closed-fields}, let $A_y := \ker((c_y)_*) < J_{\tilde{C}_y}$. Let $D_y$ be the hyperelliptic curve over $k$ obtained via (\ref{D-def}) (and Remark \ref{comp-curve}) applied to the fibers over $k$ ($D_y$ is thus the fiber at $y$ of the curve $D/S$ in Remark \ref{D/S}). The curve $D_y$ is the hyperelliptic curve of genus $g$ whose Jacobian we want to study.

Let $\tau_{(\tilde{E}_y)_\largefield}$ be a non-trivial translation which is an automorphism of $f_\largefield$, and as in Section \ref{curves-alg-closed-fields}, let $\tau_{(\tilde{C}_y)_\largefield}$ be the automorphism of $\tilde{C}_\largefield$ obtained by base-change from $\tau_{(\tilde{E}_y)_\largefield}$. Then $\largefield$ is still the smallest extension of $k$ over which $\tau_{(\tilde{C}_y)_\largefield}$ is defined.

We have $\dim(A) = [\largefield : k] = \ell -1$ and $\mathbb{Q}[(\tau^*_{\tilde{C}_\largefield})|_{A_\largefield}] \simeq \mathbb{Q}[\zeta_\ell]$. Using \cite[Proposition 2]{BDS} we conclude that $A_y$ is isogenous to the Weil restriction with respect to $\largefield|k$ of one ordinary elliptic curve over $\largefield$. It follows that  $(A_y)_\largefield$ is isogenous to a power of one ordinary elliptic curve. In particular, the Jacobian of the hyperelliptic curve $(D_y)_\largefield$, which by Lemma \ref{J_H<A} is isogenous to an abelian subvariety of $(A_y)_\largefield$, is isogenous to a power of one ordinary elliptic curve.

This proves the theorem stated in the introduction.

\section{Further examples of ordinary curves with completely decomposable Jacobian}
In the article \cite{ES} mentioned in the introduction it is not only asked whether there exist curves with completely decomposable Jacobian of arbitrarily high genus but also whether there exist such curves for every genus. To the knowledge of the authors, the question is still open in any characteristic except 2, and it is therefore interesting to give conditions on $p$ and $g$ such that a curve of genus $g$ over $\overline{\mathbb{F}}_p$ with completely decomposable Jacobian exists (for characteristic 2, see \cite{vdGvdV}). In this section, we outline how the proof of our theorem can be generalized to give further examples of ordinary curves with completely decomposable Jacobian.

By the proof of the theorem, ordinary curves with completely decomposable Jacobian exist under the following assumptions for arbitrary odd characteristic:
\begin{itemize}
\item
$g= \ell +1$,  $\ell$ an odd prime
\item
$g= (\ell-1)/2$,  $\ell$ an odd prime
\item
$g = (\ell+1)/2$,  $\ell$ an odd prime
\end{itemize}
Indeed, if $\ell$ is an odd prime, then the curves $\tilde{C}_y$ in the proof of the theorem have genus $\ell+1$ and the curves $D_y$ have genus $(\ell-1)/2$. Moreover, the curves $\tilde{C}_y$ have another canonical quotient $C_y \longrightarrow D'_y$ of degree 2 such that $g(D'_y) = (\ell+1)/2$. In the notation of Remark \ref{D/S}, the curves $D'_y$ are the normalizations of the fiber product of $\rho'_y : E'_y \longrightarrow \mathbf{P}_y$ and $u_y : \tilde{\mathbf{P}}_y \longrightarrow \mathbf{P}_y$.

\smallskip

The three conditions given above can be generalized to $\ell$ a \emph{power} of an odd prime. This follows from a generalization of the proof of the theorem. An outline is as follows: For arbitrary odd $\ell$, we have a canonical decomposition
\[ J_{\tilde{C}_y} \sim J_{C_y} \times \prod_{d | \ell, d \neq 1} (A_y)_d \; , \]
where $(A_y)_d$ is an abelian variety of dimension $\phi(d)$ with $\mathbb{Z}[\tau_{\tilde{C}_y}^*|_{(A_y)_d}] \simeq \mathbb{Z}[\zeta_d]$. Now, if $\ell$ is a prime power, the cover $Y_1 \longrightarrow Y_0$ is still cyclic, and if one chooses a point $y \in Y_0$ which is inert under this map and for which $\tilde{C}_y$ is smooth and ordinary, one can conclude with \cite[Proposition 2]{BDS} that the each of the $(A_y)_d$ is isogenous to a Weil restriction of one ordinary elliptic curve.

\smallskip

So far, we have not used the real multiplication on $J_{(D_y)_{\overline{\mathbb{F}}_p}}$ explicitly. If one does so, one can prove that there exist examples of suitable curves of genus $(\ell-1)/2$ and $(\ell+1)/2$ provided that $\ell$ is odd and $(\mathbb{Z}/\ell\mathbb{Z})^*/\langle -1 \rangle$ is \emph{cyclic} (which implies that $\ell$ is an odd prime power or the product of two odd prime powers).

We outline the proof for the genus being $(\ell-1)/2$. The case $(\ell+1)/2$ is similar.

Let $\ell$ be odd. We have a canonical decomposition
\[ J_{D_y} \sim \prod_{d | \ell, d \neq 1} (B_y)_d \; , \]
where $(B_y)_d$ is an abelian variety of dimension $\phi(d)/2$ with $\mathbb{Z}[\tau^*|_{(B_y)_d}] \simeq \mathbb{Z}[\zeta_d + \zeta_d^{-1}]$, where $\tau$ is defined as in Section \ref{curves-alg-closed-fields}.

Let us assume that $(\mathbb{Z}/\ell\mathbb{Z})^*/\langle -1 \rangle$ is cyclic. Let $\sigma_2$ be the automorphism of $Y_1 \longrightarrow Y_0$ corresponding to $-1 \in (\mathbb{Z}/\ell\mathbb{Z})^*$, and let $Y_2 := Y_1/ \langle \sigma_2 \rangle$.

Then for all $y \in Y_0$, the automorphism group of $Y_2 \longrightarrow Y_0$ operates canonically on all $\mathbb{Q}[\tau^*|_{(B_y)_d}]$. If $y$ is inert, it generates the automorphism groups of all these groups. If one chooses an inert $y$ such that $\tilde{C}_y$ (and thus $D_y$) is irreducible and ordinary, then again by \cite[Proposition 2]{BDS} each $(B_y)_d$ is isogenous to a Weil restriction of an ordinary elliptic curve.

\medskip

\subsubsection*{Acknowledgment}

With great pleasure, we thank Irene Bouw and Everett Howe for discussions and comments.
\bigskip

\bibliography{highrank-literatur}

\bibliographystyle{plain}

\bigskip

\noindent
Claus Diem, Universit\"at Leipzig, Fakult\"at f\"ur Mathematik und Informatik,
 Augustusplatz 10/11, 04109 Leipzig, Germany. diem@math.uni-leipzig.de

\smallskip

\noindent
Jasper Scholten, ESAT/COSIC, K.U. Leuven, Kasteelpark Arenberg 10, 3001 Leuven-Heverlee, Belgium. jasper.scholten@esat.kuleuven.be

\end{document}